\documentclass[11pt]{article}
\usepackage{amsfonts}
\usepackage{mathrsfs}

\usepackage[latin1]{inputenc}
\usepackage{amsmath,amssymb}
\usepackage{latexsym}
\usepackage[active]{srcltx}
\usepackage[
bookmarks=true,         
bookmarksnumbered=true, 
colorlinks=true, pdfstartview=FitV, linkcolor=blue, citecolor=blue,
urlcolor=blue]{hyperref}

\newtheorem{theorem}{Theorem}[section]
\newtheorem{corollary}[theorem]{Corollary}
\newtheorem{lemma}[theorem]{Lemma}
\newtheorem{proposition}[theorem]{Proposition}

\newtheorem{remark}[theorem]{Remark}
\numberwithin{equation}{section}
\newtheorem{example}[theorem]{Example}
\def\qed{{\hfill $\square$ \bigskip}}
\def\R{{\mathbb R}}
\def\E{{{\mathbb E}\,}}
\def\P{{\mathbb P}}

\def\Z{{\mathbb Z}}

\def\qed{{\hfill $\square$ \bigskip}}

\def\square{{\vcenter{\vbox{\hrule height.3pt
        \hbox{\vrule width.3pt height5pt \kern5pt
           \vrule width.3pt}
        \hrule height.3pt}}}}

\def\tlint{{- \kern-0.85em \int \kern-0.2em}}
\def\dlint{{- \kern-1.05em \int \kern-0.4em}}

\def\bN {{\mathbb N}}

\def \eref#1{\hbox{(\ref{#1})}}

\def \eref#1{\hbox{(\ref{#1})}}

\def\nn{{\nonumber}}

\def\nn{{\nonumber}}
\def \eref#1{\hbox{(\ref{#1})}}

\def \eref#1{\hbox{(\ref{#1})}}

\begin{document}

\title{Regularity of harmonic functions for some Markov chains with unbounded range}
\date{\empty}

\author{Fangjun Xu  \thanks{The author is supported in part by the Robert Adams Fund.}\\
\\
Department of Mathematics \\
University of Kansas \\
Lawrence, Kansas, 66045 USA}

\maketitle

\begin{abstract}
\noindent We consider a class of continuous time Markov chains on $\Z^d$. These chains are the discrete space analogue of Markov processes with jumps. Under some conditions, we show
that harmonic functions associated with these Markov chains are H\"{o}lder continuous.


\end{abstract}

\section{Introduction}

It is well known that once Harnack inequalities for Makrov processes hold, the H\"{o}lder regularity of harmonic 
functions associated with these processes follows. The technique is standard and was first developed by J. Moser in his famous paper \cite{moser}.  Recent papers \cite{bass_chen} and \cite{xu} showed that, for some singular 
 Markov processes, the H\"{o}lder regularity of harmonic functions still holds while Harnack inequalities fail. To 
 some extent, this means that Harnack inequalities are not necessary needed when proving the H\"{o}lder regularity 
 of harmonic functions. So it is natural to ask under what conditions the H\"{o}lder regularity of harmonic functions still
 holds.  In this paper, we consider a class of symmetric Markov chains defined from Dirichlet forms and then give conditions for the H\"{o}lder regularity 
 of harmonic functions associated with these Markov chains, which are the discrete space analogue of Markov processes with jumps. Our main theorem is, roughly, that an upper bound on the rate of decay of the conductances similar to that of stable processes of index $\alpha$ plus a Poincar\'{e} inequality implies that harmonic functions are H\"{o}lder continuous. We do not need a lower bound on the rate of decay of the conductances. The main difficulty here is to get near diagonal lower bounds 
 for transition densities. To obtain these lower bounds we use a scaling technique and some weighted 
 Poincar\'{e} inequalities. Scaling techniques for Markov chains and Markov processes are widely used when studying heat kernel estimates. For example, \cite{stroock_zheng}, \cite{bass_kumagai_uemura}, \cite{bass_kumagai}, \cite{chen_kumagai} and \cite{xu}. Weighted Poincar\'{e} inequalities are especially helpful when obtaining lower bounds for transition densities. See \cite{jerison}, \cite{saloff_coste} and 
 references therein. 
 
For each $x\in\Z^d$ and $A\subset\Z^d$, we define $\mu_x=1$ and
$\mu(A)=\sum\limits_{y\in A}\mu_y$. For $x\in\Z^d$ and $r>0$, let
$B(x,r)$ be the open ball in $\Z^d$ centered at $x$ with radius $r$ and $B[x,r]$ the open cube in $\Z^d$ centered at $x$ with side length $2r$. For each $x,y\in\Z^d$, let $C(x,y)$ be the conductance between $x$
and $y$. Throughout this paper, we let $\alpha\in(0,2]$ and assume that the conductance function $C(\cdot,\cdot)$ satisfies the following conditions:
\begin{enumerate}
\item[(A1)]  For any $x,y\in\Z^d$, $C(x,y)=C(y,x)\geq 0$ and $C(x,x)=0$.

\item[(A2)] There exists a positive constant $\kappa_1$ such that
\begin{equation*}
v_x=\sum_{y\in\Z^d}C(x,y)\geq\kappa_1,\quad\text{for all}\;
x\in\Z^d.
\end{equation*}

\item[(A3)] There exist positive constants $\kappa_2$ and $\kappa_3$ and a nonnegative function $\varphi:\bN\to\R+$ such that
\begin{equation*}
C(x,y)\leq\varphi(|y-x|),\; \sum_{|z|\geq r}\varphi(|z|)\leq
\frac{\kappa_2}{r^{\alpha}}\;\; \text{and}\; \sum_{|z|<r}|z|^2\varphi(|z|)\leq
\kappa_3r^{2-\alpha}
\end{equation*}
for all $x,y,z\in\Z^d$ and $r>0$.

\item[(A4)] For any open cube $B$ in $\mathbb{Z}^d$ with side length $2r$, there exist positive constants $\kappa_4$ and $\kappa_5\geq 1$ independent of $B$ such that
\begin{equation*}
\sum_B(f(x)-f_B)^2\leq \kappa_4 r^{\alpha}\sum_{\kappa_5B}\sum_{\kappa_5B}(f(y)-f(x))^2C(x,y),
\end{equation*}
where $f_B=|B|^{-1}\sum\limits_{B}f(z)$ with $|B|$ being the cardinality of $B$, and $k_5B$ is the cube with
the same center as $B$ but side length $k_5$ times as large.
\end{enumerate}

Now we use Dirichlet form to define the Markov chain associated with the conductance function $C(\cdot,\cdot)$. For each $f\in L^2(\Z^d,\mu)$, define
\begin{equation*}
\left\{
\begin{array}{rl}
\mathcal{E}(f,f)&=\; \frac{1}{2}\sum_{x,y\in\Z^d}(f(y)-f(x))^2C(x,y),\\
\\
\mathcal{F}&=\; \big\{f\in L^2(\Z^d,\mu): \mathcal{E}(f,f)<\infty\big\}.
\end{array}\right.
\end{equation*}
It is easy to see that $(\mathcal{E},\mathcal{F})$ is a regular Dirichlet form on $L^2(\Z^d,\mu)$. Let $X$ be the continuous time Markov chain corresponding to the
regular Dirichlet form $(\mathcal{E},\mathcal{F})$. In this paper, we consider Markov chains $X$ and show that the harmonic functions associated with $X$ are H\"{o}lder continuous under the assumptions (A1)-(A4).

\begin{remark} The first two conditions are mild and enable us  to define symmetric Markov chains through Dirichlet forms. The last two are used to obtain heat kernel estimates for the Markov chains. In particular, the last condition seems to be necessary for the H\"{o}lder regularity of harmonic functions associated with these Markov chains.
\end{remark}

There are some related papers, see \cite{bass_kumagai}, \cite{husseini_kassmann} and the references therein, in which the regularity of harmonic functions for Markov chains was studied. However, our results are not covered by these works. The differences between this work and  \cite{bass_kumagai}, \cite{husseini_kassmann} are given below.
\begin{itemize}
\item (A2) of \cite{bass_kumagai} implies our assumption (A4) with $\alpha=2$ and $\kappa_5\geq1$ through a comparison with the simple random walk. The conductance function $C_{xy}$ in \cite{bass_kumagai} satisfies the above assumptions (A1)-(A4) with $\alpha=2$. When $\alpha=2$, our assumption (A3) corresponds to the uniform second moment condition, which is substantive in \cite{bass_kumagai}. When $0<\alpha<2$, our assumption (A3) says that the uniform second moment condition is not needed. Even in the case $\alpha=2$, our method is a little different from that of \cite{bass_kumagai}. Bass and Kumagai in \cite{bass_kumagai} used the global weighted Poincar\'{e} inequality to obtain the near diagonal lower bound while we use local weighted Poincar\'{e} inequalities.

\item  In \cite{husseini_kassmann}, Husseini and
Kassmann considered Markov chains which are similar to stable processes. The essential assumption in \cite{husseini_kassmann} is (A3) which concerns the lower bound of the conductances. Our results do not need such assumption. See Example \ref{ex2} for conductances that do not satisfy the assumption (A3) in \cite{husseini_kassmann}. 
\end{itemize}

The paper is organized as follows.  In section 2 we obtain heat kernel estimates for $X$ and then give near diagonal lower bounds for the transition densities of $X$. In section 3, we prove a support theorem. In section 4, we show the H\"{o}lder regularity of harmonic functions associated with $X$. In section 5, we give a few examples in which assumptions (A1)-(A4) are satisfied. 

Throughout this paper, the letter c with or without a subscript indicates a positive constant whose exact value is unimportant and may change from line to line.

\bigskip

\section{Heat Kernel Estimates}
 
We start this section with the following Nash inequality.
\begin{proposition}  \label{prop2.1}
There exists $c_1$ such that
\begin{equation*}
||f||^{2+\frac{2\alpha}{d}}_2\leq
c_1\mathcal{E}(f,f)||f||^{\frac{2\alpha}{d}}_1,\quad\text{for
all}\; f\in\mathcal{F}.
\end{equation*}
\end{proposition}
Proof: For any $s>0$, let $\{Q_i\}^{\infty}_{i=1}$ be a sequence of open cubes in $\Z^d$ which have equal side length $2s$ and satisfy
\begin{equation*}
(1)~ Q_i\cap Q_j=\emptyset~\text{for}~i\neq j \quad \text{and} \quad (2)~ \cup^{\infty}_{i=1}
2Q_i=\Z^d.
\end{equation*}

From assumption (A4),
\begin{align*}
\sum_{\Z^d}f^2&\leq \sum^{\infty}_{i=1}\sum_{2Q_i}f^2\\
&\leq 2\sum^{\infty}_{i=1}\bigg(\sum_{2Q_i}(f-f_{2Q_i})^2+|2Q_i|f^2_{2Q_i}\bigg)\\
&\leq 2\sum^{\infty}_{i=1}\sum_{2Q_i}(f-f_{2Q_i})^2+2\sum^{\infty}_{i=1}|2Q_i|f^2_{2Q_i}\\
&\leq c_2s^{\alpha}\sum^{\infty}_{i=1}\sum_{2\kappa_5Q_i}\sum_{2\kappa_5Q_i}\big(f(x)-f(y)\big)^2C(x,y)+c_3s^{-d}\sum^{\infty}_{i=1}\big(\sum_{2Q_i}|f|\big)^2\\
&\leq c_4s^{\alpha}\mathcal{E}(f,f)+c_5s^{-d}\|f\|_1^2,
\end{align*}
where $f_{2Q_i}=\frac{1}{|2Q_i|}\sum\limits_{2Q_i}|f|$. Therefore, for all
$s>0$, we have
\begin{equation}
\|f\|_2^2\leq c_4s^{\alpha}\mathcal{E}(f,f)+c_5s^{-d}\|f\|_1^2.   \label{e.2.1}
\end{equation}
Choosing $s$ to minimize the right-hand side of \eref{e.2.1} completes the proof. \qed

Write $p(t,x,y)$ for the transition density of $X_t$.
\begin{proposition} \label{prop2.2} There exists $c_1$ such
that
\begin{equation*}
p(t,x,y)\leq c_1(t^{-d/\alpha}\wedge 1),\quad\text{for all}\; t>0.
\end{equation*}
\end{proposition}
Proof: It is obvious that $p(t,x,y)\leq 1$ for all $x,y\in\Z^d$ and $t>0$. Using Theorem 2.1 in
\cite{carlen_kusuoka_stroock} and Proposition \ref{prop2.1}, we know that there exists $c$ such that $p(t,x,y)\leq c t^{-d/\alpha}$ for all $x,y\in\Z^d$ and $t>0$. Combining these estimates gives the desired result. \qed

For each $\rho\geq 1$, set $\mathcal{S}=\rho^{-1}\mathbb{Z}^d$. For each $x\in\mathcal{S}$ and $A\subset\mathcal{S}$, let $\mu^{\rho}_x=\rho^{-d}$
and $\mu^{\rho}(A)=\sum\limits_{y\in A}\mu^{\rho}_y$. Define the rescaled process $V$ as
\begin{equation*}
V_t=\rho^{-1}X_{\rho^{\alpha}t},\quad t\geq 0.
\end{equation*}
Using similar arguments as in \cite{xu}, we see that the Dirichlet form corresponding to $V$ is
\begin{equation*}
\left\{
\begin{array}{rl}
\mathcal{E}^{\rho}(f,f)&=\; \sum\limits_{\mathcal{S}}\sum\limits_{\mathcal{S}}\big(f(y)-f(x)\big)^2C^{\rho}(x,y),\\
\mathcal{F}_{\rho}&=\; \big\{f\in L^2(\mathcal{S},\mu^{\rho}):\;\mathcal{E}^{\rho}(f,f)<\infty\big\},
\end{array}\right.
\end{equation*}
where $C^{\rho}(x,y)=\rho^{\alpha-d}C(\rho x,\rho y)$ for all $x,y\in\mathcal{S}$. 

Write $p^{\rho}(t,\cdot,\cdot)$ for the transition density of $V_t$. Then we have
\begin{equation} \label{e.2.2}
p^{\rho}(t,x,y)=\rho^dp(\rho^{\alpha}t,\rho x,\rho y)
\end{equation}
for all $x,y\in\mathcal{S}$ and $t\geq 0$.  The process $V$ satisfies the following Poincar\'{e} inequality. 
\begin{lemma} \label{lema2.3}
For any open cube $B$ in $\mathcal{S}$ with side length $2r$, there is a constant $c$ independent of $B$ and $\rho$ such that
\begin{equation*}
\sum_B(f(x)-f_B)^2\rho^{-d}\leq cr^{\alpha}\sum_{\kappa_5B}\sum_{\kappa_5B}(f(y)-f(x))^2C^{\rho}(x,y).
\end{equation*}
\end{lemma}
Proof: This follows from assumption (A4) and change of variables.
\qed

For $\lambda$ large enough, let $V^\lambda$ be the process $V$ with
jumps larger than $\lambda$ removed. Write
$p^{\rho,\lambda}(t,x,y)$ for the transition density of $V^{\lambda}_t$.
\begin{proposition}   \label{prop2.4}
There exists $c_1$ independent of $\rho$ and $\lambda$ such that
\begin{equation*}
p^{\rho, \lambda}(t,x,y) \leq c_1 t^{-d/\alpha} e^t .
\end{equation*}
\end{proposition}
Proof: Under the first two parts of assumption (A3), the above upper bound follows easily from Theorem 2.1 in \cite{carlen_kusuoka_stroock}, Lemma \ref{lema2.3} and the proof of Proposition \ref{prop2.1}. \qed

We can obtain a better upper bound for the transition density of $V^{\lambda}$.
\begin{lemma} \label{lema2.5}There exist $c_1$ and $c_2$ independent of $\rho$ and $\lambda$ such that
\[
p^{\rho,\lambda}(t,x,y)\leq
c_1t^{-\frac{d}{\alpha}}e^{c_2t}e^{-|x-y|/\lambda}.
\]
\end{lemma}
Proof: Applying Theorem 3.25 in \cite{carlen_kusuoka_stroock} and
Proposition \ref{prop2.4}, we have
\begin{equation}  \label{e.2.3}
p^{\rho,\lambda}(t,x,y)\leq c_2t^{-d/\alpha}e^{t-E(2t,x,y)},
\end{equation}
where
\begin{align*}
E(t,x,y)&=\sup\big\{|\psi(x)-\psi(y)|-t\Lambda(\psi)^2:\Lambda(\psi)<\infty\big\},\\
\Lambda(\psi)^2&=\|e^{-2\psi}\Gamma_\lambda[e^{\psi}]\|_{\infty}\vee\|e^{2\psi}\Gamma_\lambda[e^{-\psi}]\|_{\infty},\\
\Gamma_\lambda[v](\xi)&=\sum_{\eta\in\mathcal{S},
|\eta-\xi|\leq\lambda}\big(v(\eta)-v(\xi)\big)^2C(\rho\xi,\rho\eta)\rho^{\alpha}.
\end{align*}

Let $\psi(\xi)=\lambda^{-1}\big(|\xi-x|\wedge|y-x|\big)$. Then
$|\psi(\eta)-\psi(\xi)|\leq|\eta-\xi|/\lambda$ and
\[
\big(e^{\psi(\eta)-\psi(\xi)}-1\big)^2\leq|\psi(\eta)-\psi(\xi)|^2
e^{2|\psi(\eta)-\psi(\xi)|}\leq c_3|\eta-\xi|^2/\lambda^2
\]
for all $\eta,\xi\in \mathcal{S}$ with $|\eta-\xi|\leq\lambda$. Therefore
\[
e^{-2\psi(\xi)}\Gamma_{\lambda}[e^{\psi}](\xi)
=\sum_{\eta\in \mathcal{S}, |\xi-\eta|\leq\lambda}\big(e^{\psi(\eta)-\psi(\xi)}-1\big)^2C(\rho\eta,\rho\xi)\rho^{\alpha}\leq c_4\lambda^{-\alpha}\leq c_4.
\]
In the last second inequality we used the last part of assumption (A3). The same upper bound is obtained if $\psi$ is replaced by $-\psi$. Note
that $|\psi(x)-\psi(y)|=|x-y|/\lambda$. Substituting these estimates
into \eref{e.2.3}, we have our result after doing some algebra. \qed

For any set $A\subset\Z^d$, let
\begin{equation*}
T_A=\inf\big\{t\geq 0: X_t\notin A\big\}\quad\text{and}\quad
\tau_A=\inf\big\{t\geq 0: X_t\in A\big\}.
\end{equation*}

The upper bound in Lemma \ref{lema2.5} implies the following key exit time
estimates for $X$. The proof is the same as the
one given in Proposition 3.4 of \cite{bass_kumagai} except some minor modifications.
\begin{theorem} \label{thm2.6}
For $a>0$ and $0<b<1$, there exists $\gamma=\gamma(a,b)\in(0,1)$
such that for every $R>0$ and $x\in\Z^d$,
\begin{eqnarray*}
\P^x\big(\tau_{B(x,aR)}(X)<\gamma R^\alpha\big)\leq b.
\end{eqnarray*}
\end{theorem}

Next we are going to obtain near diagonal lower bounds for the transition densities of $X$. 

\begin{proposition} \label{prop2.7}The following two statements are equivalent:
\begin{enumerate}
\item[(1)] There is an $\epsilon$ such that
\begin{equation*}
p(t,x,y)\geq \epsilon t^{-d/\alpha}
\end{equation*}
for all $t\geq 1$ and $|y-x|\leq 2t^{1/\alpha}$.

\item[(2)] There is an $\epsilon$ such that
\begin{equation*}
p^{\rho}(1,\rho^{-1}x,\rho^{-1}y)\geq \epsilon
\end{equation*}
for all $\rho\geq 1$ and $|\rho^{-1}y-\rho^{-1}x|\leq 2$.
\end{enumerate}
\end{proposition}
Proof: This follows easily from \eref{e.2.2} and change of variables.
\qed
\begin{remark} In fact, statements (1) and (2) are also equivalent to
the following one: There is an $\epsilon$ such that
\begin{equation*}
p^{\rho}(t,\rho^{-1}x,\rho^{-1}y)\geq \epsilon t^{-d/\alpha}
\end{equation*}
for all $t\geq \rho^{-\alpha}$, $\rho\geq 1$ and $|\rho^{-1} y-\rho^{-1} x|\leq 2 t^{1/\alpha}$. 
\end{remark}

In the remainder of this section, we first prove the statement (2) in Proposition \ref{prop2.7} and then obtain the near diagonal lower bound for the transition densities of $X$.

For any $R>0$ and $x_0\in\mathcal{S}$, let $B=B[x_0,R]$ be the open cube in $\mathcal{S}$ centered at $x_0$ with side length $2R$,
\begin{equation}
\phi_R(x)=c_1\Big(R^2-|x_0-x|_m^2\Big)^{+}
\end{equation}
where $|x_0-x|_m=\max\{|x^1_0-x^1|,\cdots, |x^d_0-x^d|\}$ and $c_1$ is chosen so that $\sum\limits_{B}\phi_R(x)=\rho^d$, and set 
\begin{equation*}
\overline{f}=\sum_{B} f(x)\phi_{R}(x)\rho^{-d}.
\end{equation*}
Then we have the following local weighted Poincar\'{e} inequality with its proof given in Appendix Two.
\begin{proposition} \label{prop2.9}For any $\rho\geq 1$, there exists a constant $c_1$ independent of $\rho$ and $R$ such that
\[
\sum_{B} (f(x)-\overline{f})^2\phi_R(x)\rho^{-d}\leq c_1 R^{\alpha}\sum_{B} \sum_{B} (f(x)-f(y))^2\big(\phi_R(x)\wedge\phi_R(y)\big) C^{\rho}(x, y)
\]
for all $R\in\bigcup\limits^{\infty}_{n=0}[\frac{n}{\rho}+\frac{1}{4\rho},\frac{n}{\rho}+\frac{1}{\rho}]$.
\end{proposition}
\bigskip

We now consider $V$ killed on exiting $B$. Since
\begin{equation*}
\P^x(V_t\in A, \tau_B>t)\leq \P^x(V_t\in A)=\sum_{A}p^{\rho}(t,x,y)\mu^{\rho}_y,
\end{equation*}
this means that $\P^x(V_t=y, \tau_B>t)$ has a density bounded by
$p^{\rho}(t,x,y)$. Write $p^{\rho}_B(t,x,y)$ for the density of $\P^x(V_t=y, \tau_B>t)$.  Then we can use Proposition \ref{prop2.9} to get lower bound
for the transition density $p^{\rho}_B(1,x,y)$ when $x$ and $y$ are not far
away. See the following proposition for details.  The proof of the following proposition is long and similar to that of Proposition 4.9 in \cite{barlow_bass_chen_kassmann}, Theorem 3.4 in \cite{chen_kim_kumagai}, and Theorem 2.5 in \cite{foondun}. We postpone it to Appendix One.
\begin{proposition} \label{prop2.10} For $R\in[2d,4d]$, there exists $c_1$ independent of $\rho$, $x_0$ and $R$ such that
$$
p^{\rho}_B(1,x,y)\geq c_1,
$$
for every $(x,y)\in B(x_0,3R/4)\times B(x_0, 3R/4)$.
\end{proposition}

\begin{theorem} \label{thm2.11} There is an $\epsilon$ such that
\begin{equation*}
p(t,x,y)\geq \epsilon t^{-d/\alpha}
\end{equation*}
for all $t\geq 1$ and $|y-x|\leq 2t^{1/\alpha}$.
\end{theorem}
Proof: From the argument before Proposition \ref{prop2.10}, we see that $p^{\rho}(1,x,y)\geq p^{\rho}_B(1,x,y)$ for all $x,y\in \mathcal{S}$. Then using Propositions \ref{prop2.7} and \ref{prop2.10} gives the desired near diagonal lower bound for $p(t,x,y)$.
\qed

\bigskip

\section{Support Theorem}

\begin{lemma} \label{lema3.1} Given $\delta>0$, there exists $\kappa$ such that if
$x,y\in\Z^d$ and $A\subset\Z^d$ with dist$(x,A)$ and dist$(y,A)$
both larger than $\kappa t^{1/\alpha}$, then
\begin{equation}
\P^x\big(X_t=y, T_A\leq t\big)\leq\delta t^{-d/\alpha}.
\end{equation}
\end{lemma}
Proof:   Let
$S_A=\sup\{s\leq t: X_s\in A\}$ be the last hitting time of $A$ before
time $t$. Then
\begin{align*}
\P^x\big(X_t=y, t/2\leq T_A\leq t\big)
&\leq \P^x\big(X_t=y, t/2\leq S_A\leq t\big)\\
&= \P^y\big(X_t=x, T_A\leq t/2\big).
\end{align*}
The last equation follows from time reversal, see Lemma 4.5 of \cite{bass_kumagai}. 
Using strong Markov property and Proposition \ref{prop2.2}, we have
\begin{align*}
\P^y\big(X_t=x, T_A\leq t/2\big)&=\P^y\big(1_{\{T_A\leq
t/2\}}\P^{X_{T_A}}(X_{t-T_A}=x)\big)\nn\\
&\leq c_1 (t/2)^{-d/\alpha}\P^y\big(T_A\leq t/2\big)\nn\\
&\leq c_1 (t/2)^{-d/\alpha}\P^y(\tau_{B(y,\kappa t^{1/\alpha})}\leq
t/2)\\
&\leq\delta t^{-d/\alpha}.
\end{align*}
Here we used Theorem \ref{thm2.6} in the last inequality by choosing proper $\kappa$. Similarly,
\begin{equation*}
\P^x\big(X_t=y, T_A\leq t/2\big)\leq \delta t^{-d/\alpha}.
\end{equation*}
Combining these estimates gives our result.
\qed

\begin{proposition} \label{prop3.2} For all $t\geq 1$, there exist $c_1$ and $\theta\in(0,1)$ such
that if $|x-z|$, $|y-z|\leq t^{1/\alpha}$, $x,y,z\in\Z^d$ and
$r\geq t^{1/\alpha}/\theta$, then
\begin{equation}  \label{e.3.2}
\P^x\big(X_t=y,\tau_{B(z,r)}>t\big)\geq c_1t^{-d/\alpha}.
\end{equation}
\end{proposition}
Proof: Choose $\delta=\epsilon/2$ in Lemma \ref{lema3.1}. Then for $r>(\kappa+1)t^{1/\alpha}$ we have
\begin{align*}
&\P^x\big(X_t=y,\tau_{B(z,r)}>t\big)\\
= &\P^x\big(X_t=y\big)-\P^x\big(X_t=y,\tau_{B(z,r)}\leq t\big)\\
\geq &\frac{\epsilon}{2} t^{-d/\alpha}.
\end{align*}
Here we used Theorem \ref{thm2.11} in the last inequality. 
\qed

\begin{remark}  \label{rem3.3}
The above proposition still holds if we replace `` $|x-z|$, $|y-z|\leq t^{1/\alpha}$, $x,y,z\in\Z^d$ '' with `` $|x-y|\leq 2t^{1/\alpha}$, $x,y\in\Z^d$ '' and `` z '' in \eref{e.3.2} with `` x '', respectively.
\end{remark}

\begin{corollary} \label{cor3.4} For each $\epsilon\in(0,1)$, there exists
$\theta=\theta(\epsilon)\in(0,1)$ with the following property: if
$x,y\in\Z^d$ with $|x-y|<t^{1/\alpha}$, $t\in[1,\theta^{\alpha}
r^{\alpha})$, and $\Gamma\subset B(y,t^{1/\alpha})$ satisfies
$\mu(\Gamma)t^{-d/\alpha}\geq\epsilon$, then
\begin{equation}
\P^x(X_t\in\Gamma\;\text{and}\;\tau_{B(x,r)}>t)\geq c_1\epsilon.
\end{equation}
\end{corollary}
Proof: This follows easily from Proposition \ref{prop3.2} and Remark \ref{rem3.3}.
\qed

\begin{remark} \label{rem3.5} In fact, the condition `` $t\in[1,\theta^{\alpha}
r^{\alpha})$ '' in the above corollary can be relaxed to `` $t\in[0,\theta^{\alpha}
r^{\alpha})$ ''. 

\end{remark}

\begin{proposition}  \label{prop3.6}
For each $\epsilon\in(0,1)$, there exist constants $c_1$ and $\eta=\eta(\epsilon)\in(0,1)$
such that for any $x\in\Z^d$, if $A\subset B(x,\eta r)$ satisfies
$\mu(A)/\mu(B(x,\eta r))\geq\epsilon$, then
\begin{equation}
\P^x(T_A<\tau_{B(x,r)})\geq c_1\epsilon.
\end{equation}
\end{proposition}
Proof:  Choose $\eta=2^{-\alpha}\theta$ and $t=(\eta r)^{\alpha}$. The above proposition follows from Corollary \ref{cor3.4} and Remark \ref{rem3.5}.
\qed

\bigskip

\section{H\"{o}lder Continuity}
The following lemma can be easily proved by using Propositions \ref{prop2.2}
and \ref{prop3.2}. We refer to Lemma 5.2 in \cite{bass_kumagai} for its proof.
\begin{lemma} \label{lema4.1} There exist constants $c_1$ and $c_2$ such that
\[
c_1r^{\alpha}\leq \E^x\tau_{B(x,r)}\leq c_2r^{\alpha}.
\]
\end{lemma}

Since $X$ is a Hunt process, there is a L\'{e}vy system formula for
it. We refer to \cite{chen_kumagai} for its proof.
\begin{lemma} \label{lema4.2} For any nonnegative function $f$ on $\Z^d\times\Z^d$
that vanishes on the diagonal and any stopping time $T$,
\begin{equation*}
\E^x\Big[\sum_{s\leq
T}f(X_{s-},X_s)\Big]=\E^x\Big[\int^T_0\sum_{y\in\Z^d}f(X_s,y)C(x,y)ds\Big].
\end{equation*}
\end{lemma}

We say that $h$ is harmonic with respect to $X$ in a domain $D$ if
$h(X_{t\wedge\tau_D})$ is a $\P^x$-martingale for every $x$ in $D$.

\begin{theorem} Suppose that $h$ is bounded on $\Z^d$ and harmonic in $B(x_0,r)$
with respect to the process $X$. Then there exist constants $c$ and
$\beta\in(0,\alpha)$ such that
\begin{equation*}
|h(x)-h(y)|\leq c\bigg(\frac{|x-y|}{r}\bigg)^{\beta}\sup|h|.
\end{equation*}
\end{theorem}
Proof:  Without loss of generality, we assume that $0\leq h\leq 1$
on $\Z^d$. From Proposition \ref{prop3.6} we know that there exist constants $c_1$ and $\eta$ such that if $A\subset B(x,\eta
r)$ with $|A|/|B(x,\eta r)|\geq 1/4$, then
\begin{equation*}
\P^x(T_A<\tau_{B(x,r)})\geq c_1.
\end{equation*}

From Lemmas \ref{lema4.1} and \ref{lema4.2}, there exists $c_2$ such that
\begin{equation*}
\P^x\big(X_{\tau_{B(x,r)}}\notin B(x,s)\big)\leq
c_2\big(\frac{r}{s}\big)^{\alpha},\qquad\text{for all}\; s\geq 2r.
\end{equation*}
Let $\gamma=1-\frac{c_1}{4}$ and
$\rho=\eta\wedge\big(\frac{\gamma}{2}\big)^{1/\alpha}\wedge\Big(\frac{c_1\gamma^2}{8c_2}\Big)^{1/\alpha}$.
We need to show
\begin{equation*}
\sup_{B(x,\rho^k r)}h-\inf_{B(x,\rho^k
r)}h\leq\gamma^k,\qquad\text{for all}\; k.
\end{equation*}

For simplicity of notation, set
\begin{equation*}
B_i=B(x,\rho^i r),\;\; \tau_i=\tau_{B_i},\;\;
a_i=\sup_{B_i}h,\;\; \text{and}\;\; b_i=\inf_{B_i}h.
\end{equation*}

By the assumption that $0<h<1$ on $\Z^d$, we see $a_i-b_i\leq 1\leq\gamma^i$
for $i\leq 0$. Suppose $a_i-b_i\leq\gamma^i$ for $i\leq k$.
Now we only need to prove
\begin{equation*}
a_{k+1}-b_{k+1}\leq\gamma^{k+1}.
\end{equation*}

Notice that $b_k\leq h\leq a_k$ on $B_{k+1}$. Define
\begin{equation*}
A=\big\{z\in B(x,\rho^{k+1}r): h(z)\leq\frac{a_k+b_k}{2}\big\}.
\end{equation*}
We can assume $\mu(A)/\mu\big(B(x,\rho^{k+1}r)\big)\geq 1/2$. Otherwise we
use $1-h$ instead of $h$ in the above definition of $A$. By the
definition of $a_{k+1}$ and $b_{k+1}$, we can choose $z_1, z_2\in
B_{k+1}$ such that $a_{k+1}=h(z_1)$ and $b_{k+1}=h(z_2)$. By
optional stopping,
\begin{align*}
h(z_1)-h(z_2)
&= \E^{z_1}\big[h(X_{T_A\wedge\tau_k})-h(z_2)\big]\\
&= \E^{z_1}\big[h(X_{T_A})-h(z_2); T_A<\tau_k\big]\\
& +\E^{z_1}\big[h(X_{\tau_k})-h(z_2); T_A>\tau_k,
X_{\tau_k}\in B_{k-1}\big]\\
& +\sum^{\infty}_{i=1}\E^{z_1}\big[h(X_{\tau_k})-h(z_2);
T_A>\tau_k, X_{\tau_k}\in B_{k-1-i}-B_{k-i}\big]\\
& \leq \big(\frac{a_k+b_k}{2}-b_k\big)\E^{z_1}\big(T_A<\tau_k\big)+(a_k-b_k)\P^{z_1}\big(T_A>\tau_k\big)\\
& +\sum^{\infty}_{i=1}\big(a_{k-1-i}-b_{k-1-i}\big)\P^{z_1}\big(X_{\tau_k}\notin
B_{k-i}\big)\\
&\leq \big(a_k-b_k\big)\Big(1-\frac{\P^{z_1}(T_A<\tau_k)}{2}\Big)+\sum^{\infty}_{i=1}c_2\gamma^{k-1}\big(\rho^{\alpha}/\gamma\big)^i\\
&\leq \big(1-\frac{c_1}{2}\big)\gamma^k+2c_2\gamma^{k-2}\rho^{\alpha}\\
&\leq \big(1-\frac{c_1}{2}\big)\gamma^k+\frac{c_1}{4}\gamma^k\\
&= \gamma^{k+1}.
\end{align*}

For any $x,y\in B(x_0,r)$, let $k$ be the smallest integer such that
$|y-x|<\rho^k r$. Then $\log\big(|x-y|\big)\geq(k+1)\log\rho+\log r$ and
\begin{equation*}
\big|h(y)-h(x)\big|\leq e^{k\log\gamma}\leq
c_3e^{\big(\frac{\log\gamma}{\log\rho}\big)\log\big(\frac{|x-y|}{r}\big)}=c_3\Big(\frac{|x-y|}{r}\Big)^{\frac{\log\gamma}{\log\rho}}.
\end{equation*}
By the definition of $\gamma$ and $\rho$, it is easy to see that
$\log\gamma/\log\rho\in(0,\alpha)$. Our result follows with
$\beta=\log\gamma/\log\rho$. 
\qed

\bigskip

\section{Examples}
In this section, we give conductance functions which satisfy assumptions (A1)-(A4).
\begin{example} 
For $\alpha\in(0,2)$ and $d\geq 2$, we define the conductance functions $C_{\alpha,1}(\cdot,\cdot)$ by
\begin{equation*}
C_{\alpha,1}(x,y)=\left\{
\begin{array}{rl}
\frac{c(x,y)}{|x-y|^{d+\alpha}} &~\text{if}~ y\neq x;\\
0 &~\text{otherwise},
\end{array}\right.
\end{equation*}
where $c(x,y)=c(y,x)$ and $0<c_1\leq c(x,y)\leq c_2<\infty$ for all $x,y\in\Z^d$. The parabolic Harnack inequality holds for the Markov chains corresponding to $C_{\alpha,1}(\cdot,\cdot)$ (see, \cite{bass_Levin}). 
\end{example}

\begin{example}  \label{ex2}
For $\alpha\in(0,2)$ and $d\geq 2$, let $\mathbb{Z}_i$ be the i-th coordinate
axis in $\mathbb{Z}^d$. We define the conductance functions $C_{\alpha,2}(\cdot,\cdot)$ by
\begin{equation*}
C_{\alpha,2}(x,y)=\left\{
\begin{array}{rl}
\frac{c(x,y)}{|x-y|^{1+\alpha}} &~\text{if}~ y-x\in
\bigcup\limits^d_{i=1}\mathbb{Z}_i\backslash\{0\};\\
0 &~\text{otherwise},
\end{array}\right.
\end{equation*}
where $c(x,y)=c(y,x)$ and $0<c_1\leq c(x,y)\leq c_2<\infty$ for all $x,y\in\Z^d$. The Markov chains corresponding to $C_{\alpha,2}(\cdot,\cdot)$ are the discrete space analogue of the singular stable-like processes in \cite{bass_chen} and \cite{xu}. When $ c(x,y)\equiv1$, the Markov chain corresponding to $C_{\alpha,2}(\cdot,\cdot)$ is the discrete space analogue of the d-dimensional L\'{e}vy process whose coordinate processes are independent 1-dimensional symmetric $\alpha$-stable processes.
\end{example}

\begin{example} For $d\geq 3$, let $e^i$ be the unit vector in $\R^d$ with the $i$-th coordinate being $1$. Let $b_n=n^{n^n}$ and $a_n$ be two sequences of positive numbers with $\sum\limits^{\infty}_{n=1} a_n\leq 1/8$ and $\sum\limits^{\infty}_{n=1} a_nb^2_n< \infty$. Let $\epsilon=2\sum\limits^{\infty}_{n=1} a_n$. We define the conductance function $C_{2,3}(\cdot,\cdot)$ by
\[
C_{2,3}(x,y)=\left\{
\begin{array}{rl}
a_n &~\text{if}~ y-x=\pm b_n e^1;\\
\frac{1-\epsilon}{2(d-1)} &~\text{if}~y-x=\pm e^j\; \text{and}\; j=2,\dots,d;\\
0 &~\text{otherwise}.
\end{array}\right.
\]
This example is from \cite{bass_kumagai}. The conductance function $C_{2,3}(\cdot,\cdot)$ satisfies assumptions (A1)-(A4) with $\alpha=2$. The uniform Harnack inequality does not hold for the Markov chain corresponding to $C_{2,3}(\cdot,\cdot)$ (see, \cite{bass_kumagai}).
\end{example}

\bigskip

\section{Appendix One} 
The goal of this section is to prove Proposition \ref{prop2.10}. Recall the definition of the transition density $p^{\rho}_B(t,x,y)$ of $V$ killed upon exiting the open cube $B$ in $\mathcal{S}$ with center $x_0$ and side length $2R\in[4d,8d]$. Notice that the Dirichlet form for $V^B$ ($V$ killed upon exiting the open cube $B$) is $(\mathcal{E}^{\rho},\mathcal{F}^B_{\rho})$ where 
\begin{equation*}
\mathcal{F}^B_{\rho}=\big\{f:\mathcal{F}_{\rho}:\; f=0\; \text{on}\; B^c\big\}.
\end{equation*}
So for $f\in\mathcal{F}^B_{\rho}$,
\begin{equation*}
\mathcal{E}^{\rho}(f,f)=\sum_B\sum_B(f(x)-f(y))^2C^{\rho}(x,y)+\sum_Bf(x)^2\kappa_B(x)\mu^{\rho}_x,
\end{equation*}
where $\kappa_B(x)=2\sum\limits_{B^c}C(\rho x,\rho y)\mu^{\rho}_y$.

\begin{lemma} \label{lema6.1} There exists a positive constant $c_1$ independent of $\rho$ and $B$ such that
\begin{equation*}
p^{\rho}_B(t,x,y)\leq c_1t^{-d/\alpha}\quad\text{and}\quad \big|\frac{\partial p^{\rho}_B(t,x,y)}{\partial t}\big|\leq c_1t^{-1-\frac{d}{\alpha}}
\end{equation*}
for all $x,y\in B$ and $t>0$.
\end{lemma}
Proof: The first inequality follows immediately from Proposition \ref{prop2.2} and the argument before Proposition \ref{prop2.10}. Since
\begin{equation*}
\sum_B\sum_B p^{\rho}_B(t,x,y)^2\mu^{\rho}_x\mu^{\rho}_y \leq \sum_B p^{\rho}_B(2t,x,x)^2\mu^{\rho}_x<\infty,
\end{equation*}
the symmetric semigroup $P^B_t$ of $V^B$ is a Hilbert-Schmidt operator on $L^2(B,\mu^{\rho})$ and so it is compact 
and has a discrete spectrum $\big\{e^{-\lambda _i t}, 1\leq i\leq N \big\}$, with repetitions according to multiplicity. Here $N$ is a natural number determined by 
the Hilbert space $L^2(B,\mu^{\rho})$. Let $\big\{\phi_i, 1\leq i\leq N \big\}$ be the corresponding eigenfunctions normalized to have unit $L^2$-norm on $B$ and to be orthogonal to each other. Then
\begin{equation*}
p^{\rho}_B(t,x,y)=\sum^N_{i=1}e^{-\lambda_i t}\phi_i(x)\phi_i(y).
\end{equation*}
Hence $p^{\rho}_B(t,x,y)$ is differential with respect to $t$ and 
\begin{align*}
|\frac{\partial p^{\rho}_B(t,x,y)}{\partial t}|
&= |-\sum^N_{i=1}\lambda_i e^{-\lambda_i t}\phi_i(x)\phi_i(y)|\\
&\leq \sum^N_{i=1}\lambda_ie^{-\lambda_i t/2} e^{-\lambda_i t/2}|\phi_i(x)| |\phi_i(y)|\\
&\leq \frac{c_2}{t}\Big(\sum^N_{i=1}e^{-\lambda_i t/2}\phi_i(x)^2\Big)^{1/2}\Big(\sum^N_{i=1}e^{-\lambda_i t/2}\phi_i(y)^2\Big)^{1/2}\\
&= \frac{c_2}{t}\Big(p^{\rho}_B(t/2,x,x)\Big)^{1/2}\Big(p^{\rho}_B(t/2,y,y)\Big)^{1/2}\\
&\leq c_3t^{-1-\frac{d}{\alpha}}.
\end{align*}
Here we used the fact that $h(x)=xe^{-xt/2}$ is bounded on $[0,\infty)$ by $c_2/t$.

For $\epsilon\in(0,1)$, define 
\begin{equation*}
G(t)=\sum\phi_R(x)\log p^{\rho,\epsilon}_B(t,x,y_0)\mu^{\rho}_x,
\end{equation*}
where $p^{\rho,\epsilon}_B(t,x,y)=p^{\rho}_B(t,x,y)+\epsilon$.

\begin{lemma} Fix $y_0\in B$. Then, for every $t>0$,
\begin{equation*}
G'(t)=-\mathcal{E}^{\rho}\Big(p^{\rho}_B(t,\cdot,y_0),\frac{\phi_R(\cdot)}{p^{\rho,\epsilon}_B(t,\cdot,y_0)}\Big).
\end{equation*}
\end{lemma}
Proof: From Lemma 1.3.3 of \cite{fukushima_oshida_takeda} and Lemma \ref{lema6.1}, we see that $p^{\rho}_B(t,x,y_0)$ as a function of $x\in B$ is in $\mathcal{F}^B_{\rho}$.  By Lemma 1.3.4 of \cite{fukushima_oshida_takeda}, we have
\begin{eqnarray*}
&&-\mathcal{E}^{\rho}\Big(p^{\rho}_B(t,\cdot,y_0),\frac{\phi_R(\cdot)}{p^{\rho,\epsilon}_B(t,\cdot,y_0)}\Big)\\
&=&\lim_{h\downarrow 0}\frac{1}{h}\Big(p^{\rho}_B(t+h,\cdot,y_0)-p^{\rho}_B(t,\cdot,y_0),\frac{\phi_R(\cdot)}{p^{\rho,\epsilon}_B(t,\cdot,y_0)}\Big)\\
&=&\lim_{h\downarrow 0}\frac{1}{h}\Big(p^{\rho,\epsilon}_B(t+h,\cdot,y_0)-p^{\rho,\epsilon}_B(t,\cdot,y_0),\frac{\phi_R(\cdot)}{p^{\rho,\epsilon}_B(t,\cdot,y_0)}\Big)\\
&=&\lim_{h\downarrow 0}\frac{1}{h}\sum\Big(\frac{p^{\rho,\epsilon}_B(t+h,x,y_0)}{p^{\rho,\epsilon}_B(t,x,y_0)}-1\Big)\phi_R(x) \mu^{\rho}_x.
\end{eqnarray*}

Moreover, 
\begin{equation*}
G'(t)=\lim_{h\to 0}\frac{1}{h}\sum\Big(\log p^{\rho,\epsilon}_B(t+h,x,y_0)-\log p^{\rho,\epsilon}_B(t,x,y_0) \Big)\phi_R(x) \mu^{\rho}_x.
\end{equation*}
Let 
\begin{equation*}
F(h)=\Big[\log p^{\rho,\epsilon}_B(t+h,x,y_0)-\log p^{\rho,\epsilon}_B(t,x,y_0)-\Big(\frac{p^{\rho,\epsilon}_B(t+h,x,y_0)}{p^{\rho,\epsilon}_B(t,x,y_0)}-1\Big)\Big]\phi_R(x) \mu^{\rho}_x.
\end{equation*}
Then 
\begin{equation*}
F'(h)=\frac{\partial p^{\rho,\epsilon}_B(t,x,y_0)}{\partial t}\Big(p^{\rho,\epsilon}_B(t,x,y_0)-p^{\rho,\epsilon}_B(t+h,x,y_0)\Big)\frac{\phi_R(x) }{p^{\rho,\epsilon}_B(t+h,x,y_0)p^{\rho,\epsilon}_B(t,x,y_0)}.
\end{equation*}
Now the lemma follows easily from using the mean value theorem, Lemma \ref{lema6.1} and the dominated convergence theorem.
\qed
\newline
{\bf Proof of Proposition 2.10}: Recall that $R\in[2d,4d]$. With the help of the above results and Proposition \ref{prop2.9}, Proposition \ref{prop2.10} follows from similar arguments as in Proposition 4.9 of \cite{barlow_bass_chen_kassmann}, Theorem 3.4 of \cite{chen_kim_kumagai}, or Theorem 2.5 of \cite{foondun}. 
\qed

\bigskip

\section{Appendix Two} 

In this appendix, we prove Proposition \ref{prop2.9}. If $B$ is an open cube in $\mathcal{S}$, we define $\overline{B}$ to be the union of all closed cubes in $\R^d$ with centers in $B$ and equal side length $\rho^{-1}$, and $\tilde{B}$ to be the interior of $\overline{B}$. If $f$ is defined on $\mathcal{S}$, we define $\tilde{f}$ as the extension of $f$ to $\R^d$:
\begin{equation*}
\tilde{f}(x)=f([x]_{\rho}),
\end{equation*}
where $[x]_{\rho}=(\rho^{-1}[\rho x^1], \dots, \rho^{-1}[\rho x^d])$ for $x=(x^1, \dots, x^d)\in\R^d$. Similarly, we can define $\tilde{C}^{\rho}(\cdot,\cdot)$ as the extension of $C^{\rho}(\cdot,\cdot)$ to $\R^d\times\R^d$.

With the above notation, the Poincar\'{e} inequality in Lemma \ref{lema2.3} can be written as follows.
\begin{lemma} \label{lema7.1} For any open cube $B$ in $\mathcal{S}$ with side length $2r$, there is a constant $c$ independent of $\rho$ and $B$ such that 
\begin{equation*}
\int_{\tilde{B}}(\tilde{f}(x)-\tilde{f}_{\tilde B})^2\, dx \leq c\, r^{\alpha}\int_{\widetilde{\kappa_5B}}\int_{\widetilde{\kappa_5B}}(\tilde{f}(y)-\tilde{f}(x))^2\tilde{C}^{\rho}(x,y)\rho^{2d}\, dx\, dy,
\end{equation*}
where $\tilde{f}_{\tilde B}=|\tilde B|^{-1}\int_{\tilde B} \tilde{f}(z)\, dz$ and $|\tilde B|$ is the Lebesgue measure of $\tilde B$.
\end{lemma}

Then we have the following result.
\begin{lemma} For any open cube $B$ in $\R^d$ with side length $2r$, there is a constant $c$ independent of $\rho$ and $B$ such that 
\begin{equation*}
\int_{B}(\tilde{f}(x)-\tilde{f}_{B})^2\, dx \leq c\, r^{\alpha}\int_{2B}\int_{2B}(\tilde{f}(y)-\tilde{f}(x))^2\tilde{C}^{\rho} (x,y)\rho^{2d}\, dx\, dy.
\end{equation*}
\end{lemma}
Proof: Lemma \ref{lema7.1} implies that there are constants $c_1$ and $k>1$ such that 
\begin{equation*}
\int_{B}(\tilde{f}(x)-\tilde{f}_{B})^2\, dx\leq c_1r^{\alpha}\int_{kB}\int_{kB}(\tilde{f}(y)-\tilde{f}(x))^2\tilde{C}^{\rho}(x,y)\rho^{2d}\, dx\, dy.
\end{equation*}
Our result then follows from the Jerison's technique in \cite{jerison} or a well-known argument mentioned in \S5.3.1 of \cite{saloff_coste}.
\qed

Using similar arguments as in the proof of Theorem 5.3.4 in \cite{saloff_coste}, we obtain the following weighted Poincar\'{e} inequality.
\begin{proposition} \label{prop7.3} For any open cube $B$ in $\R^d$ with center in $\mathcal{S}$ and side length $2R$, there is a constant $c$ independent of $
R$ and $\rho$ such that
\begin{equation*} 
\int_{B} (\tilde{f}(x)-\overline{\tilde{f}})^2\phi_R(x)dx \leq c R^{\alpha}\int_{B} \int_{B} (\tilde{f}(x)-\tilde{f}(y))^2\phi_R(x)\wedge\phi_R(y) \tilde{C}^{\rho} (x, y)\rho^{2d}\, dx\, dy,
\end{equation*}
where $\overline{\tilde{f}}=\int_B\tilde{f}(x)\phi_R(x)dx$.
\end{proposition}

{\bf Proof of Proposition 2.9}: We see that Proposition \ref{prop2.9} is trivial when $R\leq1/\rho$ since both sides of the inequality equal zero. For all $R$ with $R-[R]\in[\frac{1}{4\rho},\frac{1}{2\rho}]$, by Proposition \ref{prop7.3}, we obtain that 
\begin{align*}
\sum_{B} (f(x)-\overline{f})^2\phi_R(x) \rho^{-d}
&\leq\sum_{B} (f(x)-\overline{\tilde{f}})^2\phi_R(x)\rho^{-d}\\
&\leq2^d\int_{B} (\tilde{f}(x)-\overline{\tilde{f}})^2\phi_R(x)dx\\
&\leq c_1 R^{\alpha}\int_{B} \int_{B} (\tilde{f}(x)-\tilde{f}(y))^2\phi_R(x)\wedge\phi_R(y) \tilde{C}^{\rho}(x, y)\rho^{2d}\, dx\, dy\\
&\leq c_2 R^{\alpha}\sum_{B} \sum_{B} (f(x)-f(y))^2\phi_R(x)\wedge\phi_R(y) C^{\rho}(x, y),
\end{align*}
where the sets $B$ in the second and third inequalities are open cubes in $\R^d$ instead of $\mathcal{S}$. This implies that there exists a constant $c_3$ independent of $\rho$ and $R$ such that
\begin{equation*}
\sum_{B} (f(x)-\overline{f})^2\phi_R(x) \rho^{-d}\leq c_3R^{\alpha}\sum_{B} \sum_{B} (f(x)-f(y))^2\phi_R(x)\wedge\phi_R(y) C^{\rho}(x, y)
\end{equation*}
for all $R\in\bigcup\limits^{\infty}_{n=0}[\frac{n}{\rho}+\frac{1}{4\rho}, \frac{n}{\rho}+\frac{1}{2\rho}]$. It is easy to see that the above inequality also holds when $R\in\bigcup\limits^{\infty}_{n=0}[\frac{n}{\rho}+\frac{1}{2\rho},\frac{n}{\rho}+\frac{1}{\rho}]$. \qed

\section*{Acknowledgements} 
The author would like to thank Professor Richard Bass for suggesting this problem. The author also thanks anonymous referees and the associate editor for helpful comments.

\end{document}